\documentclass{amsart}
\usepackage[curve,matrix,arrow]{xy}
\usepackage{amssymb, amscd}
\newtheorem{theorem}{Theorem}
\newtheorem{lemma}[theorem]{Lemma}

\newtheorem{proposition}[theorem]{Proposition}

\newtheorem{remark}[theorem]{Remark}

\newtheorem{numbered equation}[theorem]{(\hspace*{-0.1cm}}

\DeclareMathOperator{\Ev}{Ev_0}

\DeclareMathOperator{\co}{C_0}

\newcommand{\opPsi}{\overline{\Theta}}

\newcommand{\bQ}{{\mathbb Q}}
\newcommand{\bQt}{\widetilde{\mathbb Q}}

\newcommand{\bP}{{\mathbb P}}
\newcommand{\hh}{{\mathbb H}}
\newcommand{\mZ}{{\mathbb Z}}
\newcommand{\bZ}{{\mathbb Z}}

\newcommand{\D}{{\mathcal D}}

\newcommand{\spec}{Sp^{\Sigma}}

\newcommand{\Ch}{{\mathcal C}h}

\renewcommand{\to}{\longrightarrow}

\newcommand{\varrow}[1]{\hbox to #1{\rightarrowfill}}
\newcommand{\varl}[2]{\stackrel{#2}{\hbox to #1{\leftarrowfill}}}

\newcommand{\varrx}[2]{\stackrel{#2}{\hbox to #1{\rightarrowfill}}}
\newcommand{\parallelarrows}[1]{\begin{array}{c} {\hbox to
#1{\rightarrowfill}}  \vspace{-0.35cm} \\ {\hbox to
#1{\rightarrowfill}} \end{array}}

\begin{document}

\title[Correction for ``DGAs and $H{\mathbb Z}$-algebras"] 
{Correction to ``$H{\mZ}$-algebra spectra are differential graded algebras" 
}

\date{\today; 2000 AMS Math.\ Subj.\ Class.: 55P43, 18G35, 55P62, 18D10} 
\author{Brooke Shipley}
\thanks{Research partially supported by NSF Grants No. 0417206 and No. 0706877}
\address{Department of Mathematics \\ 508 SEO m/c 249\\   
851 S. Morgan Street\\ Chicago, IL 60607-7045 \\ USA}
\email{bshipley@math.uic.edu}

\begin{abstract} 
This correction article is actually unnecessary.
The proof of 
Theorem 1.2, concerning commutative $H\bQ$-algebra
spectra and commutative differential graded algebras, in the author's 
paper [{\em American
Journal of Mathematics} {\bf 129} (2007) 351-379 (arxiv:math/0209215v4)] 
is correct as originally stated.
Neil Strickland carefully proved that $D$ is 
symmetric monoidal~\cite{neil};
so Proposition 4.7  and hence also Theorem 1.2 hold as stated.   
Strickland's proof will appear in joint work 
with Stefan Schwede~\cite{neil-stefan}; see related work in 
[arxiv:0810.1747]~\cite{neil-arxiv}.   Note here $D$ is defined as a colimit
of chain complexes; in contrast, non-symmetric monoidal functors analogous to 
$D$ are defined as homotopy colimits of spaces in previous work of the 
author~\cite{thh-shipley}.

We leave the old alternate approach to Theorem 1.2 below, 
with expository changes in the introduction, 
since it does provide another slightly weaker, non-natural statement.      

\end{abstract}  
\maketitle

In the author's paper~\cite{S}, the proof of Theorem 1.2 
is correct as stated; the functor $D$ 
is symmetric monoidal. 
The author's confusion about this fact came from the 
comparison of this functor $D$, which is defined as a colimit of chain
complexes, with the functor $D$ in~\cite{thh-shipley} which
is defined as a homotopy colimit of spaces.  
See also the discussion of commutative $I$-monoids in section 2.2
of~\cite{CS}.  In the topological case $D$ is not symmetric monoidal;
in the algebraic case though the functor $D$ is symmetric. 

Since the paper~\cite{S} is mainly concerned with associative algebras, 
the only place this issue arises is in the proof
of Theorem 1.2.  
As stated in Remark 2.11 in~\cite{S}, the main theorems (Theorem 1.1,
Corollary 2.15 and Corollary 2.16 in~\cite{S}) would also hold with the ``three
step" functors $H$ and $\Theta$ replaced by the ``four step"
functors $\overline{\hh} = U L c \co f F_0 c$
and $\opPsi  = \Ev f i \phi^* N Z c $
where $c$ and $f$ are the appropriate cofibrant and fibrant replacement
functors.  Since here the functors $\Ev$, $i$, $\phi^* N$ and $Z$ are
symmetric monoidal, we have the following non-natural version
of Theorem 1.2 from~\cite{S} with $\Theta$ replaced by $\opPsi$. This 
statement first appeared as Theorem 1.3 in~\cite{S-early}.  

\begin{theorem}\label{thm-Q-comm}
For $C$ any commutative $H\bQ$-algebra, $\opPsi C$ is weakly equivalent
to a commutative differential graded $\bQ$-algebra.
\end{theorem}

\begin{proof}
As noted in the proof of Theorem 1.2 from~\cite{S}, the reason
$\opPsi$ is not symmetric monoidal is because the cofibrant and
fibrant replacement functors involved in $\opPsi$ are not
symmetric monoidal.   This is why $\opPsi C$ is only weakly
equivalent and not isomorphic to a commutative dg $\bQ$ algebra.

The method for dealing with the cofibrant replacement functor
in $\opPsi$ proceeds as in~\cite{S}.  As proved there, a natural
zig-zag of weak equivalences exists between $Z c$ and the
symmetric monoidal functor $\alpha^* \bQt$.
Let $\opPsi'  = \Ev f i \phi^* N \alpha^* \bQt $.  Then $\opPsi C$
is naturally weakly equivalent to $\opPsi' C$.  

Next we need to consider the fibrant replacement functor $f$ which
appears in $\opPsi'$ (and $\opPsi$).   This is the fibrant replacement
functor in the model category of monoids in $\spec(\Ch_{\bQ})$.   
As in~\cite{S-early}, we exchange $f$ for the fibrant
replacement functor $f'$ in the model category of commutative monoids
in $\spec(\Ch_{\bQ})$ as established below in Proposition~\ref{prop-com}.
For any commutative monoid $A$ in $\spec(\Ch_{\bQ})$, we thus have
two weak equivalences $A \to fA$ and $A \to f'A$.   Since $f'A$ is also fibrant
as a monoid and $A \to fA$ is a trivial cofibration of monoids, 
lifting provides a weak equivalence $fA \to f'A$.  If we let $\opPsi'' C = 
\Ev f' i \phi^* N \alpha^* \bQt $,  we then have a (non-natural) weak
equivalence $\opPsi'C \to \opPsi'' C$.  Since $\opPsi C$
is weakly equivalent to $\opPsi' C$ and $\opPsi'' C$ is 
a commutative differential graded $\bQ$-algebra, this completes the proof. 
\end{proof}

\begin{remark}
{\em Although Theorem~\ref{thm-Q-comm} does not give a natural
identification of $\opPsi C$ with a commutative DGA, for any
small, fixed $I$-diagram ${\D}$ of commutative $H{\bQ}$-algebra spectra
there will be a map of $I$-diagrams from $\opPsi \D$ to  
an $I$-diagram of commutative DGAs which is given by a variant 
of $\opPsi'' \D$ with $f'$ replaced by the fibrant replacement
functor in the model category of $I$-diagrams of commutative monoids
in $\spec(\Ch_{\bQ})$ given by~\cite[11.6.1]{hh}.  
}\end{remark}


\begin{proposition}\label{prop-com}
There is a model category structure on the category of commutative
monoids in $\spec(\Ch_{\bQ})$ in which a map is a weak equivalence 
or fibration if and only if the underlying map in $\spec(\Ch_{\bQ})$
is so. 
\end{proposition}

Let $S_{\bQ}$ denote the unit and let $\otimes_S$ denote the monoidal product 
in $\spec(\Ch_{\bQ})$. 
To establish this model category we use the lifting property 
from~\cite[2.3(i)]{ss00} applied to the {\em free commutative monoid}
functor $\bP$ which is left adjoint to the forgetful functor
from commutative monoids in $\spec(\Ch_{\bQ})$ to the underlying
object in $\spec(\Ch_{\bQ})$.  Namely, $\bP(M)
= \bigvee_{n\geq 0} M^{(n)}/\Sigma_n$ where $M^{(n)} = M \otimes_S
\cdots \otimes_S M$ is the $n$th
tensor power of $M$ over $S_{\bQ}$.

Let $I$ denote the generating cofibrations and $J$ denote the
generating trivial cofibrations in $\spec(\Ch_{\bQ})$;  
see~\cite[7]{hovey}.   
To establish the lifting criterion 
in~\cite[2.3]{ss00}, we first show that applying $\bP$ to any
map in $J$ produces a stable equivalence. 
We do this by showing that in the source and target the orbit
constructions can be replaced by homotopy orbits without
changing the homotopy type.  

\begin{lemma}\label{lem-es}
Let $X$, $Y$  be in $\spec(\Ch_{\bQ})$ and $n \geq 1$.
\begin{enumerate}
\item The map 
\[ E \Sigma_n \otimes_{\Sigma_n} X^{(n)}\to 
  X^{(n)}/\Sigma_n  \]
is a level equivalence.
\item The map
\[ (E \Sigma_n \otimes_{\Sigma_n} X^{(n)}) \otimes_S Y \to
 ( X^{(n)}/\Sigma_n) \otimes_S Y \]
is also a level equivalence.
\end{enumerate}
\end{lemma}

\begin{proof}
The first statement follows directly from the fact that given any 
$\Sigma_n$-equivariant complex $A$ in $\Ch_{\bQ}$, then  
\[ E \Sigma_n \otimes_{\Sigma_n} A \to A/\Sigma_n  \]
is a quasi-isomorphism.  The second statement follows as well
by extending the $\Sigma_n$-action trivially to $Y$ and shifting
the parentheses. 
\end{proof}

Next we show that pushouts of maps in $\bP(J)$ are stable equivalences
and level cofibrations.  Since directed colimits of such maps are
again stable equivalences, Proposition~\ref{prop-com} then follows from
Lemmas~\ref{lem-es} and~\ref{lem-po} by~\cite[2.3]{ss00}. 

\begin{lemma}\label{lem-po}
Let $f:T \to U$ be a cofibration in $\spec(\Ch_{\bQ})$ and
$V$ be a $\bP T$-module.  Then 
the map $q: V \to V \otimes_{\bP T} \bP U$ is a level cofibration.  
If $f$ is a trivial cofibration, then $q$ is a stable equivalence. 
\end{lemma}

\begin{proof}
This follows from the fitration arguments of~\cite[7.5, 8.6]{mandell} using
Lemma~\ref{lem-es} instead of~\cite[8.2, 8.10]{mandell}. 
Note, here one does not need to restrict to the positive cofibrant objects
since no such restriction is needed in Lemma~\ref{lem-es}. 
\end{proof}

\end{document}